\documentclass[12pt]{article}
\usepackage{amsmath,amsfonts,amssymb}
\newtheorem{thm}{Theorem}

\newtheorem{lem}[thm]{Lemma}

\newtheorem{cor}{Corollary}

\newcommand\la{\lambda}

\def\qed{\nobreak\quad\raise -2pt\hbox{\vrule\vbox to 10pt{\hrule width 6pt
\vfill\hrule}\vrule}\par\vspace{2ex}}
\def\qel{\nobreak\quad\raise -2pt\hbox{\vrule\vbox to 10pt{\hrule width 6pt
\vfill\hrule}\vrule}}
\begin{document}

\begin{center}
{\Large\bf Applications of Waring's formula to some identities of
Chebyshev polynomials}
\end{center}

\vskip 2mm \centerline{ Jiang Zeng$^{1,2}$ and Jin Zhou$^2$}

\begin{center} \small $^1$ Institut Girard Desargues,
Universit\'e Claude Bernard (Lyon I)\\
69622 Villeurbanne Cedex, France \\
{\tt zeng@igd.univ-lyon1.fr}\\
and\\
$^2$ Center for Combinatorics, LPMC,
Nankai University\\
 Tianjin 300071, People's Republic of China\\
{\tt jinjinzhou@hotmail.com} \\
\vspace{10pt}
\small AMS Math Subject Classification Numbers: 11B39, 33C05, 05E05\\
\end{center}

 \noindent
\textbf{Abstract.} Some identities of Chebyshev polynomials are
deduced from Waring's formula on symmetric functions. In
particular, these formulae generalize some recent results of
Grabner and Prodinger.
\section{Introduction}

Given a set of variables $X=\{x_1,x_2,\ldots\}$, the $k$th ($k\geq
0$) \emph{elementary symmetric polynomial} $e_{k}(X)$ is defined
by $e_0(X)=1$,
$$
e_k(X)=\sum_{i_1<\ldots<i_k}x_{i_1}\ldots x_{i_k},\quad
\text{for}\quad k\geqslant 1,
$$
 and  the $k$th  ($k\geq
0$)\emph{power sum symmetric polynomial} $p_{k}(X)$ is defined by
$p_0(X)=1$,
$$
p_k(X)=\sum_ix_i^k,\quad \text{for}\quad k\geqslant 1.
$$
Let $\la=1^{m_1}2^{m_2}\ldots $ be a partition of $n$, i.e., $m_1
1+m_2 2+\ldots +m_n n=n$, where $m_i\geq 0$ for $i=1,2,\ldots n$.
Set $l(\la)=m_1+m_2+\ldots +m_n$. According to the
\emph{fundamental theorem of symmetric polynomials}, any symmetric
polynomial can be written uniquely as a polynomial of elementary
symmetric polynomials $e_i(X)$ ($i\geq 0$). In particular, for the
power sum $p_k(x)$, the corresponding formula  is usually
attributed to Waring~\cite{CLY,Mac} and reads as follows:
\begin{equation}\label{war}
p_k(X)=\sum_{\la}(-1)^{k-l(\lambda)}\frac{k(l(\lambda)-1)!}
{\prod_{i}{m_{i}!}}e_{1}(X)^{m_1}e_{2}(X)^{m_2}\ldots,
\end{equation}
where the sum is over all the partitions
$\la=1^{m_1}2^{m_2}\ldots$ of $k$.

In a recent paper~\cite{GP} Grabner and Prodinger proved some
identities about Chebyshev polynomials using generating functions,
the aim of this paper is to show that Waring's formula provides
a natural generalization of such kind of identities.

Let $U_n$ and $V_n$ be two
sequences defined by the following recurrence relations:
\begin{align}
U_n&=pU_{n-1}-U_{n-2},&U_0=0, U_1=1,\\
V_n&=pV_{n-1}-V_{n-2},&V_0=2, V_1=p.
\end{align}
Hence $U_n$ and $V_n$ are rescaled versions of the first and
second kind of Chebyshev polynomials ${\cal U}_n(x)$ and ${\cal
T}_n(x)$, respectively:
$$
{\cal U}_n(x)=U_{n+1}(2x),\quad {\cal T}_n(x)=\frac{1}{2}T_n(x).
$$
\begin{thm} For integers
$m,n\geq 0$, let $W_n=aU_n+bV_n$ and  $\Omega=a^2+4b^2-b^2p^2$. Then the following identity holds
\begin{equation}\label{eq:wm}
W_n^{2k}+W_{n+m}^{2k}=\sum_{r=0}^{k}\theta_{k,r}(m)\Omega^{k-r}W_n^rW_{n+m}^r,
\end{equation}
where
$$
\theta_{k,r}(m)=\sum_{0\leqslant 2j\leqslant
k}(-1)^j\frac{k(k-j-1)!}{j!(k-r)!(r-2j)!}V_m^{r-2j}U_m^{2k-2r}.
$$
\end{thm}
Note that  the identities of Grabner and Prodinger~\cite{GP} correspond to the $m=1$
and implicitly $m=2$ cases of Theorem~1 (cf. Section~3).
\section{Proof of Theorem~1}
We first check the $k=1$ case of (\ref{eq:wm}):
\begin{equation}\label{fund}
W_n^2+W_{n+m}^2=V_mW_nW_{n+m}+U_m^2\Omega.
\end{equation}
Set $\alpha=(p+\sqrt{p^2-4})/2$ and $\beta=(p-\sqrt{p^2-4})/2$
then it is easy to see that
$$
U_n=\frac{\alpha^{n}-\beta^{n}}{\alpha-\beta},\quad
V_{n}=\alpha^n+\beta^n,
$$
it follows that
\[
W_n=aU_n+bV_n=A\alpha^n+B\beta^n,
\]
where $A=b+a/(\alpha-\beta)$ and $B=b-a/(\alpha-\beta)$. Therefore
\begin{align*}
V_mW_nW_{n+m}+U_m^2\Omega
&=(\alpha^m+\beta^m)(A\alpha^n+B\beta^n)(A\alpha^{n+m}+B\beta^{n+m})\\
&+\left(\frac{\alpha^m-\beta^m}{\alpha-\beta}\right)^2(a^2+4b^2-b^2p^2),
\end{align*}
 which is readily seen to be equal to $W_n^2+W_{n+m}^2$.

Next we  take the alphabet $X=\{W_n^2, W_{n+m}^2\}$, then the left-hand
side of (\ref{eq:wm}) is the power sum $p_k(X)$. On the other
hand, since
$$
e_1(X)=W_n^2+W_{n+m}^2,\quad  e_2(X)=W_n^2W_{n+m}^2,\quad e_i(X)=0
\quad \textrm{if} \quad i\geqslant3,
$$
the summation at the right-hand side of (\ref{war}) reduces to the
partitions $\la=(1^{k-2j}\,2^j)$, with  $j\geq 0$. Now, using
(\ref{fund}) Waring's formula~(\ref{war}) infers that
\begin{align*}
W_n^{2k}&+W_{n+m}^{2k} \\
 &=\sum_{0\leqslant2 j\leqslant
k}(-1)^{j} \frac{k(k-j-1)!}{j!(k-2j)!}
(V_mW_nW_{n+m}+U_m^2\Omega)^{k-2j}(W_n^2W_{n+m}^2)^{j}\\
&=\sum_{0\leqslant2j\leqslant k}\sum_{i=0}^{k-2j}(-1)^j
\frac{k(k-j-1)!}{j!i!(k-2j-i)!}V_m^{k-2j-i}
U_m^{2i}\Omega^{i}(W_{n}W_{n+m})^{k-i}
\end{align*}
Setting $k-i=r$ and exchanging the order of summations yields
(\ref{eq:wm}). \qed

\section{Some special cases}
When $m=1$ or 2, as $U_1=1$, $V_1=p$ and $U_2=p$, $V_2=p^2-2$ the
coefficient $\theta_{k,r}(r)$ of Theorem~1 is much simpler.
\begin{cor} We have
\begin{eqnarray}
\theta_{k,r}(1)&=&\sum_{0\leqslant2j\leqslant
r}(-1)^j\frac{k(k-1-j)!}{(k-r)!j!(r-2j)!}p^{r-2j},\label{coeff1}\\
\theta_{k,r}(2) &=&\sum_{0\leqslant2j\leqslant
k}(-1)^j\frac{k(k-j-1)!}{j!(k-r)!(r-2j)!}
(p^2-2)^{r-2j}p^{2k-2r}.\label{coeff2}
\end{eqnarray}
\end{cor}

We notice  that (\ref{coeff1}) is exactly the formula given by Grabner and Prodinger~\cite{GP}
for $\theta_{k,r}(1)$, while for $\theta_{k,r}(2)$  they
give a more involved formula than  (\ref{coeff2})  as follows:
\begin{cor}[Grabner and Prodinger~\cite{GP}] There holds
\begin{equation}\label{coeffGP}
\theta_{k,r}(2) =\sum_{0\leqslant\lambda\leqslant
k}(-1)^{\lambda}p^{2k-2\lambda}
\frac{k(k-\lfloor\frac{\lambda}{2}\rfloor-1)!
2^{\lceil\frac{\lambda}{2}\rceil}}{(k-r)!\lambda!(r-\lambda)!}
\prod_{i=0}^{\lfloor\frac{\lambda}{2}\rfloor-1}
(2k-2\lceil\frac{\lambda}{2}\rceil-1-2i).
\end{equation}
\end{cor}
In order to identify (\ref{coeff2}) and (\ref{coeffGP}), we need
the following identity.
\begin{lem}\label{key} We have
\begin{eqnarray}\label{fa}
&&\sum_{i=0}^{j/2}(-1)^{i}\frac{(k-i-1)!2^{j-2i}}{(j-2i)!i!}\nonumber\\
&&\hspace{1cm}= \frac{(k-\lfloor{j/2}\rfloor-1)!}{j!}2^{\lceil
j/2\rceil} \prod_{i=0}^{\lfloor j/2\rfloor-1}(2k-2\lceil
j/2\rceil-1-2i).
\end{eqnarray}
\end{lem}
\begin{pr} For $n\geq 0$ let $(a)_n=a(a+1)\ldots (a+n-1)$, then
the Chu-Vandermonde formula~\cite[p.212]{GKP} reads:
\begin{equation}\label{cv}
{}_2F_1( -n,a; c; 1):=\sum_{k\geqslant0}\frac{(-n)_k(a)_k}{(c)_k
k!}=\frac{(c-a)_n}{(c)_n}.
\end{equation}
Note that $n!=(1)_n$, so using the simple transformation formulae:
$$
(a)_{2n}=\left(\frac{a}{2}\right)_n\left(\frac{a+1}{2}\right)_{n}2^{2n},\quad
(a)_{2n+1}=\left(\frac{a}{2}\right)_{n+1}
\left(\frac{a+1}{2}\right)_{n}2^{2n+1},
$$
and
$$
(a)_{N-n}=\frac{(a)_N}{(a+N-n)_n}=(-1)^n\frac{(a)_N}{(-a-N+1)_n},
$$
 we can rewrite the left-hand side of identity~(\ref{fa}) as follows:
 \begin{eqnarray*}
  \begin{cases}
    \frac{(k-1)!}{(\frac{1}{2})_m(1)_m}
    \,{}_2F_1(-m, -m+\frac{1}{2};-k+1; 1)& \text{if $j=2m$}, \\
    \\
    \frac{(k-1)!}{{(\frac{1}{2})}_{m+1}(1)_m}
    \,{}_2F_1(-m, -m-\frac{1}{2};-k+1; 1)& \text{if $j=2m+1$},
  \end{cases}
\end{eqnarray*}
which is clearly equal to the right-hand side of (\ref{fa}) in
view of (\ref{cv}). \qed
\end{pr}

Now, expanding  the right-hand side of (\ref{coeff2}) by binomial
formula yields
$$
\sum_{0\leqslant2j\leqslant k}(-1)^j\frac{k(k-j-1)!}
{j!(k-r)!(r-2j)!}\sum_{i=0}^{r-2j}{r-2j\choose i}
p^{2i}(-2)^{r-2j-i}p^{2k-2r}.
$$
Writing  $\la=r-i$, so $\la\leq r\leq k$, and exchanging the order
of summations, the above quantity becomes
$$
\sum_{0\leqslant \la\leqslant k} (-1)^{\la} p^{2k-2\la
}\frac{k}{(k-r)!(r-\la)!}\sum_{0\leqslant j\leqslant
k/2}(-1)^j\frac{(k-j-1)!2^{\la-2j}}{(\la-2j)!j!},
$$
which yields (\ref{coeffGP}) by applying Lemma~\ref{key}.


\begin{thebibliography}{99}
\bibitem{CLY} William Y. C. Chen, Ko-Wei Lih, and  Yeong-Nan Yeh:
\emph{Cyclic Tableaux and Symmetric Functions}, Studies in Applied
Math., 94 (1995), 327-339.
\bibitem{GKP} Ronald L. Graham, Donald E. Knuth and Oren Patashnik:
\emph{Concrete Mathematics}, Addion-Wesley Pubilshing Co. 1989.
\bibitem{GP} Peter J. Grabner and Helmut Prodinger:
\emph{Some identities for Chebyshev polynomials}, Portugalia
Mathematicae 59 (2002), 311-314.
\bibitem{Mac} P. A. MacMahon:
\emph{Combinatory analysis}, Chelsea Publishing Co. New York,
1960.

\end{thebibliography}
\end{document}